\documentclass[12pt]{article}

\newcommand{\norm}{\|}
\newtheorem{theorem}{\bf Theorem}
\newtheorem{proposition}{\bf Proposition}

\newtheorem{rmk}{\bf Remark}
\newenvironment{remark}{\begin{rmk}\rm}{\end{rmk}}
\newtheorem{ex}{\bf Exercise}

\newcommand{\be}{\begin{equation}}
\newcommand{\ee}{\end{equation}}

\newcommand{\eps}{\varepsilon}

\newcommand{\conv}{{\rm conv}}

\newcommand {\CC}{{I\kern-.6em C}}
\newcommand {\NN}{{I\kern-.3em N}}

\newcommand{\bth}{\begin{theorem}}
\newcommand{\eth}{\end{theorem}}
\newcommand{\bqn}{\begin{eqnarray*}}
\newcommand{\eqn}{\end{eqnarray*}}

\begin{document}
\title{Intrinsic Volumes of the Brownian Motion Body}
\author{{\sc Fuchang Gao} \\
{\it University of Idaho}
\and
{\sc Richard A. Vitale}\\
{\it University of Connecticut}}

\date{\today}
\maketitle

\begin{abstract}
Motivated from Gaussian processes, we derive the intrinsic volumes of the infinite--dimensional {\it Brownian motion body}. The method is by discretization to a class of orthoschemes. Numerical support is offered for a conjecture of Sangwine-Yager, and another conjecture is offered on the rate of decay of intrinsic volume sequences.   
\end{abstract}

\section{Introduction}
As is well known, intrinsic volumes of convex bodies are conveniently normalized versions of quermassintegrals and are important in many settings (\cite{sang}, \cite{sch}). Unfortunately, closed form expressions are known in only a handful of cases (\cite{had49}, \cite{had55}, \cite[pp. 220--221]{had57}, \cite[pp. 224--232]{sant}). Here we present a new example. Original motivation came from the theory of Gaussian random processes, where intrinsic volumes have been a useful tool, for instance, in using the mean width of a compact convex body in a Hilbert space to characterize the so-called $GB$ ({\it G}aussian {\it B}ounded) sets (\cite{che}, \cite{sud}, \cite{tsiI}, \cite{tsiII}, \cite{tsiIII}, \cite{vit96}, \cite{vit99}). The intrinsic volume of an infinite-dimensional convex body in Hilbert space is defined to be the corresponding supremum of intrinsic volumes over all included finite-dimensional convex bodies.

Of all Gaussian processes, Brownian motion is arguably the most important. The corresponding convex body $K_B$, called the {\it Brownian motion body}, is the closed convex hull of a so-called crinkled arc $B,$ which maps $[0,1]$ continuously into Hilbert space such that $B(0)=0$ and $\norm B(\hat{t})-B(t) \norm=\sqrt{\hat{t}-t}$ for $0\leq t <\hat{t}\leq 1$ (which entails the crinkled property: for any $0<t_1<t_2<t_3<t_4$, $[B(t_2)-B(t_1)] \perp [B(t_4)-B(t_3)]$). By isometry, any such map $B(\cdot)$ can be taken: one example has $B(t)$ the indicator function of the interval $[0,t]$ regarded as an element of $L_2[0,1];$ then $K_B$ is the set of non-increasing functions in $L_2[0,1]$ that are bounded between 0 and 1.

We will show that in the case of the Brownian motion body the sequence of intrinsic volumes has a particularly attractive form:

\begin{theorem} \label{main}For $k=1,2,\ldots$,
$$V_k(K_B)=\omega_k/k!\,\,,$$
where $\omega_k=\pi^{k/2}/\Gamma(k/2 +1)$ is the volume of the $k$-dimensional unit ball. 
\end{theorem}
\noindent{\bf Remark} Note that the intrinsic volume sequence is infinite, as is the case for general non-finite-dimensional convex bodies. We return to this point in the final section.

The argument for Theorem 1 is by discretization and, in fact, nearly all the work will be to establish the following:
\begin{theorem}
In ${\bf R}^n$, let $K$ be the simplex with vertices $(0,0,\dots, 0), (1,0, \dots,0), (1,1,0,\dots,0),\linebreak (1,1,1,0\dots,0), ..., (1,1,\dots,1)$. For $1\leq k\leq n$, 
$$V_k(K)=\frac{1}{k!}\sum_{A}\frac{1}{\sqrt{l_1\, l_2 \,l_3 \cdots l_k}}\,\, ,$$
where $A$ is the set of integral vectors $(l_1,l_2,\dots, l_k)$ with $1\leq l_i$ for $1\leq i\leq k$ and $\sum_{1\leq i\leq k}l_i\leq n$.  
\end{theorem}

In the next section, we make some preliminary comments and then exhibit direct constructions for Theorem 2 in low dimension. Section 3 carries the proof of the general case and numerical support for a conjecture of Sangwine-Yager. Section 4 gives the deduction of Theorem 1 from Theorem 2 and a conjecture on the tail behavior of intrinsic volume sequences.

\section{Preliminary Comments and Direct Constructions}
We begin by recalling some facts: for $1\leq k\leq n$, let $J=\{i_0,i_1,i_2,\dots, i_k\}$ be a set of positive integers, such that $0\leq i_0<i_1<i_2<\cdots<i_k\leq n$. Let $F_J$ be the $k$-dimensional face of $K$ that contains $P_i=(\underbrace{1,1,\ldots,1}_i,0,0,\ldots,0)$, $i\in J$, and let $A_J$ denote the $k$-volume of $F_J$. Then $A_J=\sqrt{(i_k-i_{k-1})(i_{k-1}-i_{k-2})\cdots (i_1-i_0)}/k!\,\,$. It is well-known (e.g., \cite{mcm}) that
$$V_k(K)=\sum_{J}(\mbox{mes}_k\,F_J)\,\gamma_J=\sum_J A_J\,\gamma_J,$$
where $\gamma_J$ is the Gaussian measure of the normal cone $N(F_J,K)$ at $F_J$ and the summation is over all $k$-dimensional faces of $K$. To distingush it from $\gamma_J$, we use $\Gamma_J$ to denote the surface area measure of the solid angle formed by $N(F_J,K)$.

The evaluation of high dimensional solid angles is typically difficult. For three dimensional solid angles, the following old formula (probably due to Euler) can be helpful (see e.g. \cite{Eri}).
\begin{proposition}
Let {\bf a,\,b\,,c} be the unit vectors of the three extreme rays of the solid angle $E$ in ${\bf R}^3$. Then 
$$\tan (\Gamma/2)=\frac{|{\bf a}\cdot ({\bf b}\times {\bf c})|}{1+{\bf b}\cdot {\bf c}+{\bf c}\cdot {\bf a}+{\bf a}\cdot {\bf b}},$$
where $\Gamma$ is the spherical area measure of $E$.
\end{proposition}

We proceed now as follows: for $n=3$, we find $V_1(K)$ by a direct computation; for $n=4$, we use Proposition 1 to compute the required solid angles. These computations are instructive for the more elaborate argument in general dimension.

\bigskip
\noindent $\bf n=3$
\bigskip

We have $P_0=(0,0,0)$, $P_1=(1,0,0)$, $P_2=(1,1,0)$ and $P_3=(1,1,1)$. There are three edges of length $1$, namely $P_0P_1$, $P_1P_2$ and $P_2P_3$. The (interior) angles of the simplex at those edges are respectively $\pi/4$, $\pi/2$ and $\pi/4$; the exterior angles are respectively $3\pi/4$, $\pi/2$ and $3\pi/4$, and the Gaussian measures of these exterior angles are respectively $3/8$, $1/4$ and $3/8$. There are two edges of length $\sqrt{2}$, namely $P_0P_2$ and $P_1P_3$, the corresponding exterior angles are $\pi/2$ and $\pi/2$ and the Gaussian measures $1/4$ and $1/4$. There is one edge of length $\sqrt{3}$, that is $P_0P_3$. The corresponding exterior angle is $2\pi/3$. The Gaussian measure is $1/3$. Thus, 
$$V_1(K)=\frac{3}{8}+\frac{1}{4}+\frac{3}{8}+\sqrt{2}( \frac{1}{4}+\frac{1}{4})+\sqrt{3}\cdot \frac{1}{3}=1+\frac{1}{\sqrt{2}}+\frac{1}{\sqrt{3}}.$$

\bigskip
\noindent $\bf n=4$
\bigskip

Here the simplex has four edges of length $1$, namely $P_0P_1$, $P_1P_2$, $P_2P_3$ and $P_3P_4$. Let us consider first edge $P_0P_1$. There are three $3$-dimensional faces of $K$ which contain $P_0P_1$. They are $F_{0,1,2,3}$, $F_{0,1,2,4}$ and $F_{0,1,3,4}$, where $F_{0,1,2,3}$ means the face containing the vertex $P_0$, $P_1$, $P_2$ and $P_3$. 

It is easy to check that the outward normal vectors of these faces are $(0,0,0,-1)$, \linebreak $(0,0,-1/\sqrt{2}, 1/\sqrt{2})$ and $(0, -1/\sqrt{2},1/\sqrt{2}, 0)$ respectively. Embedding 
these vectors into a three dimensional hyperplane $\{(w,x,y,z): w=0\}$ and applying Proposition 1, we obtain $\tan(\Gamma_{01}/2)=-\sqrt{2}-1,$ which implies $\Gamma_{01}=5\pi/4$. The Gaussian measure is $\gamma_{01}=5/16$. Similar computations lead to the following:
\bqn
\begin{array}{|l|llll|}
\cline{1-5}
\mbox{edge length}=1&
\gamma_{01}=5/16&
\gamma_{12}=3/16&
\gamma_{23}=3/16&
\gamma_{34}=5/16\\
\cline{1-5}
\mbox{edge length}=\sqrt{2}&
\gamma_{02}=3/16&
\gamma_{13}=1/8&
\gamma_{24}=3/16&\\
\cline{1-5}
\mbox{edge length}=\sqrt{3}&
\gamma_{03}=1/6&
\gamma_{14}=1/6&&\\
\cline{1-5}
\mbox{edge length}=2&
\gamma_{04}=1/4&&&\\
\cline{1-5}
\end{array}
\eqn
Therefore 
\bqn
\sum_{0\leq i<j\leq 4}A_{ij}\gamma_{ij}&=&1\times(\gamma_{01}+\gamma_{12}+\gamma_{23}+\gamma_{34})+\sqrt{2}\times(\gamma_{02}+\gamma_{13}+\gamma_{24})\\
&&\hspace{1cm}+\sqrt{3}\times(\gamma_{03}+\gamma_{14})+\sqrt{4}\times\gamma_{04}\\
&=&\left(\frac{5}{16}+\frac{3}{16}+\frac{3}{16}+\frac{5}{16}\right)+\sqrt{2}\left(\frac{3}{16}+\frac{1}{8}+\frac{3}{16}\right)+\sqrt{3}\left(\frac{1}{6}+\frac{1}{6}\right)+\frac{1}{2}\\
&=& 1+\frac{1}{\sqrt{2}}+\frac{1}{\sqrt{3}}+\frac{1}{\sqrt{4}}.
\eqn

\section{General Dimension}
In this section, we find intrinsic volumes of $K$ in arbitrary dimension. In evaluating $\sum_J A_J\gamma_J$, we can no longer use Proposition 1, so another approach is required. Note though that for the case $n=4$, what we really used were the sums $\gamma_{01}+\gamma_{12}+\gamma_{23}+\gamma_{34}$, $\gamma_{02}+\gamma_{13}+\gamma_{24}$, and so forth. This suggests grouping the angles and finding the sums of their measures {\it within groups}. To this end, let us first see what those angles are.

For $1<i<n$, define ${\bf u}_i$ to be the $n$-dimensional vector whose $i$-th coordinate is $-1/\sqrt{2}$, $(i+1)$-st coordinate is $1/\sqrt{2}$, and the other coordinates are zeros; let also 
${\bf u}_0= ( 1,0,0,\dots,0)$ and ${\bf u}_n=( 0,0,\dots,0,-1 )$. Fix $J=\{i_0,i_1,\dots, i_k\}$. It is not hard to check that the extreme rays of $N(F_J,K)$ are the vectors ${\bf u}_i$, $i\notin J$. A good way to present this solid angle is to put the extreme ray vectors in the following matrix form:

\smallskip
{\footnotesize
\bqn
\left[
%\begin{array}{lcc|ccc|c|ccc|ccr}
\begin{array}{lcccccccccccr}
1&&&&&&&&&&&&\\
\frac{-1}{\sqrt{2}}\,\,\frac{1}{\sqrt{2}}&&&&&&&&&&&&\\
&\!\!\!\!\!\ddots\!\!\!\!\!&&&&&&&&&&&\\
&&\frac{-1}{\sqrt{2}}\,\,\frac{1}{\sqrt{2}}&&&&&&&&&&\\
\cline{1-13}
&&&\frac{-1}{\sqrt{2}}\,\,\frac{1}{\sqrt{2}}&&&&&&&&&\\
&&&&\!\!\!\!\!\ddots\!\!\!\!\!&&&&&&&&\\
&&&&&\frac{-1}{\sqrt{2}}\,\,\frac{1}{\sqrt{2}}&&&&&&&\\
\cline{1-13}
&&&&&&\!\!\ddots\!\!&&&&&&\\
\cline{1-13}
&&&&&&&\frac{-1}{\sqrt{2}}\,\,\frac{1}{\sqrt{2}}&&&&&\\
&&&&&&&&\!\!\!\!\!\!\ddots\!\!\!\!\!&&&&\\
&&&&&&&&&\frac{-1}{\sqrt{2}}\,\,\frac{1}{\sqrt{2}}&&&\\
\cline{1-13}
&&&&&&&&&&\frac{-1}{\sqrt{2}}\,\,\frac{1}{\sqrt{2}}&&\\
&&&&&&&&&&&\!\!\!\!\!\!\ddots\!\!\!\!\!&\\
&&&&&&&&&&&&\frac{-1}{\sqrt{2}}\,\,\frac{1}{\sqrt{2}}\\
&&&&&&&&&&&&-1
\end{array}
\right]
\eqn
}where the lines in the middle indicate the missing vectors ${\bf u}_i$, $i\in J$. There are $(k+1)$ small matrices appearing on the diagonal. The $k-1$ small matrices in the middle are of the same form but possibly different order. Call them $B_1$, $B_2$, $\dots$, $B_{k-1}$. First, let us look at $B_1$. It is an $(i_1-i_0-1)\times (i_1-i_0)$ matrix. The rows are $(i_1-i_0)$-dimensional vectors. We denote them by ${\bf v}_1$, ${\bf v}_2$, $\dots$, ${\bf v}_{i_1-i_0-1}$. They can be thought of as the extreme rays of an $(i_1-i_0-1)$-dimensional solid angle. For convenience, simply call $B_1$ an $(i_1-i_0-1)$-dimensional angle. The main idea is to produce $(i_1-i_0)$ angles such that (i) each has the same measure as  $B_1$, and (ii) their union can be regarded as a partition (no common interior points) of ${\bf R}^{i_1-i_0-1}$. To see this, for $1\leq l\leq i_1-i_0-1$, we delete the $l$-th row of $B_1$ and add the $(i_1-i_0)$-dimensional vector ${\bf v}_{i_1-i_0}=( 1/\sqrt{2},0,\dots,0,-1/\sqrt{2})$ at the bottom. Then we move the first $l-1$ rows to the bottom. Call the new matrix $B_1^l$. To see that it has the same measure as a solid angle as $B_1$, we display the matrices: 
\bqn
B_1=\left[
\begin{array}{rrrrrrr}
\frac{-1}{\sqrt{2}}&\frac{1}{\sqrt{2}}&0&\dots&0&0&0\\
0&\frac{-1}{\sqrt{2}}&\frac{1}{\sqrt{2}}&\dots&0&0&0\\
&&&\ddots&&&\\
0&0&0&\dots&\frac{-1}{\sqrt{2}}&\frac{1}{\sqrt{2}}&0\\
0&0&0&\dots&0&\frac{1}{\sqrt{2}}&\frac{-1}{\sqrt{2}}
\end{array}
\right]\\
\\
B_1^1=\left[
\begin{array}{rrrrrrr}
0&\frac{-1}{\sqrt{2}}&\frac{1}{\sqrt{2}}&0&\dots&0&0\\
0&0&\frac{-1}{\sqrt{2}}&\frac{1}{\sqrt{2}}&\dots&0&0\\
&&&&\ddots&&\\
0&0&0&0&\dots&\frac{-1}{\sqrt{2}}&\frac{1}{\sqrt{2}}\\
\frac{-1}{\sqrt{2}}&0&0&0&\dots&0&\frac{1}{\sqrt{2}}
\end{array}
\right]\\
\\
B_1^2=\left[
\begin{array}{rrrrrrr}
0&0&\frac{-1}{\sqrt{2}}&\frac{1}{\sqrt{2}}&\dots&0&0\\
&&&&\ddots&&\\
0&0&0&0&\dots&\frac{-1}{\sqrt{2}}&\frac{1}{\sqrt{2}}\\
\frac{-1}{\sqrt{2}}&0&0&0&\dots&0&\frac{1}{\sqrt{2}}\\
\frac{-1}{\sqrt{2}}&\frac{1}{\sqrt{2}}&0&0&\dots&0&0
\end{array}
\right]
\eqn
$$
\cdots\cdots
$$
\bqn
B_1^{i_1-i_0-1}=\left[
\begin{array}{rrrrrrr}
\frac{1}{\sqrt{2}}&0&0&\dots&0&0&\frac{-1}{\sqrt{2}}\\
\frac{-1}{\sqrt{2}}&\frac{1}{\sqrt{2}}&0&\dots&0&0&0\\
0&\frac{-1}{\sqrt{2}}&\frac{1}{\sqrt{2}}&\dots&0&0&0\\
&&&\ddots&&&\\
0&0&0&\dots&\frac{-1}{\sqrt{2}}&\frac{1}{\sqrt{2}}&0
\end{array}
\right]\hspace{.3in}\\
\eqn
Each matrix $B_1^{l}$ is a column permutation of $B_1$. Therefore the corresponding angle is a reflection of the original angle. Thus the $(i_1-i_0)$ angles are of the same measure. To show that they form a partition of ${\bf R}^{i_1-i_0-1}$, let us first argue that they are disjoint (no common interior points). For convenience, we denote $B_1^0=B_1$. Note that, for $0\leq l<m\leq i_1-i_0-1$, $B_1^l$ and $B_1^m$ do not have common interior points. In fact, if ${\bf x}=( x_1,x_2,\cdots,x_n)$ is a vector in $B_1^l$, then it is a convex combination of ${\bf v}_1$, ${\bf v}_2$, $\dots$,${\bf v}_{l-1}$,${\bf v}_{l+1}$, $\cdots$, ${\bf v}_{i_1-i_0}$, say
 $${\bf x}=a_1{\bf v}_1+\cdots+a_{l-1}{\bf v}_{l-1}+a_{l+1}{\bf v}_{l+1}+\cdots+a_{i_1-i_0}{\bf v}_{i_1-i_0}.$$
Thus 
\bqn
x_{l+1}+x_{l+2}+\cdots+x_m&=&\left(-a_{l+1}\frac{1}{\sqrt{2}}\right)+
\left(a_{l+1}\frac{1}{\sqrt{2}}-a_{l+2}\frac{1}{\sqrt{2}}\right)+\cdots+\left(a_{m-1}\frac{1}{\sqrt{2}}-a_m\frac{1}{\sqrt{2}}\right)\\
&=&-a_m\frac{1}{\sqrt{2}}<0.
\eqn
On the other hand, if ${\bf x}$ is a vector in $B_1^m$, then it is a convex combination of ${\bf v}_1$, ${\bf v}_2$, $\dots$,${\bf v}_{m-1}$,${\bf v}_{m+1}$, $\cdots$, ${\bf v}_{i_1-i_0}$, say 
$${\bf x}=b_1{\bf v}_1+\cdots+b_{m-1}{\bf v}_{m-1}+b_{m+1}{\bf v}_{m+1}+\cdots+b_{i_1-i_0}{\bf v}_{i_1-i_0}.$$ Thus 
\bqn
x_{l+1}+x_{l+2}+\cdots+x_m&=&\left(a_l\frac{1}{\sqrt{2}}-a_{l+1}\frac{1}{\sqrt{2}}\right)+
\left(a_{l+1}\frac{1}{\sqrt{2}}-a_{l+2}\frac{1}{\sqrt{2}}\right)+\cdots+\left(a_{m-1}\frac{1}{\sqrt{2}}\right)\\
&=&a_l\frac{1}{\sqrt{2}}>0.
\eqn
But this is a contradiction. Therefore the $B_1^l$'s are disjoint. To see that they form a partition of ${\bf R}^{i_1-i_0-1}$, consider the $(i_1-i_0-1)$-dimensional hyperplane $${\bf H}^{i_1-i_0}=\{(x_1,x_2,\dots,x_{i_1-i_0}):\, x_1+x_2+\cdots +x_{i_1-i_0}=0\}.$$ 
Because the vector sequence ${\bf v}_1$, ${\bf v}_2$, $\dots$, ${\bf v}_{i_1-i_0}$ has rank $i_1-i_0-1$, the hyperplane ${\bf H}^{i_1-i_0}$ is the linear span of these vectors. For any ${\bf x}\in {\bf H}^{i_1-i_0}$, ${\bf x}$ can be written as ${\bf x}=c_1{\bf v}_1+c_2{\bf v}_{2}+\cdots+c_{i_1-i_0}{\bf v}_{i_1-i_0}.$ Without loss of generality, assume $c_1\leq c_2\leq \cdots\leq c_{i_1-i_0}$. Then
\bqn
{\bf x}&=&c_1({\bf v}_1+{\bf v}_2+\cdots+{\bf v}_{i_1-i_0})+(c_2-c_1){\bf v}_{2}+\cdots+(c_{i_1-i_0}-c_1){\bf v}_{i_1-i_0}\\
&=&(c_2-c_1){\bf v}_{2}+\cdots+(c_{i_1-i_0}-c_1){\bf v}_{i_1-i_0}.
\eqn
This means that ${\bf x}$ is in the angle $B_1^1$. Thus $B_1$, $B_1^1$, $B_1^2$, $\dots$, $B_1^{i_1-i_0-1}$ form a partition of ${\bf H}^{i_1-i_0-1}$. 

Similarly, look at block $B_p$, $1<p<k$. We can produce $i_p-i_{p-1}$ angles, each of the same measure as the angle $B_p$, such that they form a partition of ${\bf R}^{i_p-i_{p-1}-1}$. 

Now we return to the big matrix. The row vectors from different blocks are orthogonal. If we replace any small matrix $B_1$ by $B_1^l$, or $B_2$ by $B_2^m$, $\dots$, the new angle is a reflection of the original one. Following the argument above, we can find $(i_1-i_0)\times (i_2-i_1)\times \cdots \times (i_k-i_{k-1})$ angles with the same measure, such that they form a partition of $E_{i_0}^{i_k-i_0}\times {\bf R}^{i_k-i_0-k}$, where $E_{i_0}^{i_k-i_0}$ is the $(n-i_k+i_0-1)$-angle determined by the matrix 
\bqn
\left[
\begin{array}{rrrrrlrrrrrr}
1&&&&&&&&&&\\
\frac{-1}{\sqrt{2}}&\frac{1}{\sqrt{2}}&&&&&&&&&&\\
&&\ddots&&&&&&&&&\\
&&&\frac{-1}{\sqrt{2}}&\frac{1}{\sqrt{2}}&&&&&&&\\
\cline{1-12}
&&&&&&\frac{-1}{\sqrt{2}}&\frac{1}{\sqrt{2}}&&&&\\
&&&&&&&\frac{-1}{\sqrt{2}}&\frac{1}{\sqrt{2}}&&&\\
&&&&&&&&&\ddots&&\\
&&&&&&&&&&\frac{-1}{\sqrt{2}}&\frac{1}{\sqrt{2}}\\
&&&&&&&&&&&-1
\end{array}
\right].
\eqn
This matrix is obtained by deleting the middle blocks of the big matrix. For $1<i<n-i_k+i_0-1$, let ${\bf w}_i$ be the $(n-i_k+i_0)$-dimensional vector whose $i$-th coordinate is $-1/\sqrt{2}$, $(i+1)$-st coordinate is $1/\sqrt{2}$, and the other coordinates are zeros; also, let ${\bf w}_0= ( 1,0,0,\dots,0)$ and ${\bf w}_{n-i_k+i_0-1}=( 0,0,\dots,0,-1 )$. Thus the extreme rays of $E_{i_0}^{i_k-i_0}$ consist of all the vectors ${\bf w}_i$, except ${\bf w}_{i_0}$. As we commented at the beginning of the section, we only need to find the sums of the measures within groups. Here we will not focus on finding the measure of each individual $E_{i_0}^{i_k-i_0}$, instead, we show that the angles $E_0^{i_k-i_0}$, $E_1^{i_k-i_0}$, $\dots$, $E_{n-i_k+i_0-1}^{i_k-i_0}$ form a partition of ${\bf R}^{n-i_k+i_0-1}$. 

First, these angles are disjoint. Suppose ${\bf x}=( x_1,x_2,\dots,x_{n-i_k+i_0-1})$ is an interior point of both $E_{l}^{i_k-i_0}$ and $E_{m}^{i_k-i_0}$, $0\leq l<m\leq n-i_k+i_0-1$. Then ${\bf x}$ can be expressed as a convex combination 
$$
{\bf x}=a_0{\bf w}_0+\cdots+a_{l-1}{\bf w}_{l-1}+a_{l+1}{\bf w}_{l+1}+\cdots+a_{n-i_k+i_0-1}{\bf w}_{n-i_k+i_0-1}.
$$ 
If $m<n-i_k+i_0-1$, then
\bqn
x_{l+1}+x_{l+2}+\cdots+x_m&=&\left(-a_{l+1}\frac{1}{\sqrt{2}}\right)+
\left(a_{l+1}\frac{1}{\sqrt{2}}-a_{l+2}\frac{1}{\sqrt{2}}\right)+\cdots+\left(a_{m-1}\frac{1}{\sqrt{2}}-a_m\frac{1}{\sqrt{2}}\right)\\
&=&-a_{l+1}\frac{1}{\sqrt{2}}-a_m\frac{1}{\sqrt{2}}<0.
\eqn
If $m=n-i_k+i_0-1$, then 
\bqn
x_{l+1}+x_{l+2}+\cdots+x_m&=&\left(-a_{l+1}\frac{1}{\sqrt{2}}\right)+
\left(a_{l+1}\frac{1}{\sqrt{2}}-a_{l+2}\frac{1}{\sqrt{2}}\right)+\cdots+\left(a_{m-1}\frac{1}{\sqrt{2}}-a_m\right)\\
&=&-a_{l+1}\frac{1}{\sqrt{2}}-a_m<0.
\eqn
On the other hand, because ${\bf x}$ is also assumed to be an interior point of $E_{m}^{i_k-i_0}$, ${\bf x}$ can be expressed as a convex combination 
$$
{\bf x}=b_0{\bf w}_0+\cdots+b_{m-1}{\bf w}_{m-1}+b_{m+1}{\bf w}_{m+1}+\cdots+b_{n-i_k+i_0-1}{\bf w}_{n-i_k+i_0-1}.
$$ 
Thus, if $l>1$, then 
\bqn
x_{l+1}+x_{l+2}+\cdots+x_m&=&\left(a_{l}\frac{1}{\sqrt{2}}-a_{l+1}\frac{1}{\sqrt{2}}\right)+
\left(a_{l+1}\frac{1}{\sqrt{2}}-a_{l+2}\frac{1}{\sqrt{2}}\right)+\cdots+\left(a_{m-1}\frac{1}{\sqrt{2}}\right)\\
&=&a_{l}\frac{1}{\sqrt{2}}+a_{m-1}\frac{1}{\sqrt{2}}>0.
\eqn
If $l=1$, then
\bqn
x_{l+1}+x_{l+2}+\cdots+x_m&=&\left(a_{l}-a_{l+1}\frac{1}{\sqrt{2}}\right)+
\left(a_{l+1}\frac{1}{\sqrt{2}}-a_{l+2}\frac{1}{\sqrt{2}}\right)+\cdots+\left(a_{m-1}\frac{1}{\sqrt{2}}\right)\\
&=&a_{l}+a_{m-1}\frac{1}{\sqrt{2}}>0.
\eqn
In all cases, we obtain a contradiction. This means that the angles $E_0^{i_k-i_0}$, $E_1^{i_k-i_0}$, $\dots$, $E_{n-i_k+i_0-1}^{i_k-i_0}$ are disjoint. 

Now we show that those angles form a partition of ${\bf R}^{n-i_k+i_0-1}$. Because 
the vector sequence $\{{\bf w}_i\}$, $0\leq i\leq n-i_k+i_0-1$, has rank $n-i_k+i_0-1$, their linear span is ${\bf R}^{n-i_k+i_0-1}$. For any ${\bf x}\in {\bf R}^{n-i_k+i_0-1}$, ${\bf x}$ can be written as 
$$
{\bf x}=\frac{c_0}{\sqrt{2}}{\bf w_0}+c_1{\bf w}_1+c_2{\bf w}_{2}+\cdots+\frac{c_{n-i_k+i_0-1}}{\sqrt{2}}{\bf w}_{n-i_k+i_0-1}.
$$
Without loss of generality, we assume $c_0\leq c_1\leq c_2\leq \cdots\leq c_{n-i_k+i_0-1}$. Then
\bqn
{\bf x}&=&c_0\left(\frac{{\bf w}_0}{\sqrt{2}}+{\bf w}_1+{\bf w}_2+\cdots+\frac{{\bf w}_{n-i_k+i_0-1}}{\sqrt{2}}\right)\\
&& +(c_1-c_0){\bf w_{1}}+\cdots+(c_{n-i_k+i_0-1}-c_0){\bf w}_{n-i_k+i_0-1}\\
&=&(c_1-c_0){\bf w}_{1}+\cdots+(c_{n-i_k+i_0-1}-c_0){\bf w}_{n-i_k+i_0-1}.
\eqn
This means that ${\bf x}$ is a convex combination of ${\bf w}_i$, $1\leq i\leq n-i_k+i_0-1$, which implies that ${\bf x}$ is in the angle $E_0^{i_k-i_0}$. Thus we have proved that $E_0^{i_k-i_0}$, $E_1^{i_k-i_0}$, $\dots$, $E_{n-i_k+i_0-1}^{i_k-i_0}$ form a partition of ${\bf R}^{n-i_k+i_0-1}$.

If we let ${\rm Gauss}(E_{i_0}^{i_k-i_0})$ be the Gaussian measure of $E_{i_0}^{i_k-i_0}$, then the angle and angle measure identifications we have shown above imply that
$$\gamma_J=\frac{1}{(i_1-i_0)\times (i_2-i_1)\times \cdots \times (i_k-i_{k-1})}\cdot {\rm Gauss}(E_{i_0}^{i_k-i_0}),$$
and for any fixed $i_k-i_0$,
$$\sum_{i=0}^{n-i_k+i_0-1}{\rm Gauss}(E_i^{i_k-i_0})=1.$$
To conclude the proof of Theorem 2, it remains to observe that 
\bqn
\sum_{J}A_J\gamma_J&=&\sum_{J}\frac{\sqrt{(i_1-i_0)\times (i_2-i_1)\times \cdots \times (i_k-i_{k-1})}}{k!}\gamma_J\\
&=&\sum_{l_1+l_2+\cdots+l_k\leq n}\frac{1}{k!}\cdot \frac{1}{\sqrt{l_1\,l_2\,\cdots \times l_k}}.
\eqn

\begin{remark}
Sangwine-Yager \cite{sang} has conjectured that if $K$ is a convex set in ${\bf R}^n$ and if $a_1\leq a_2\leq \cdots\leq a_n$ are the real parts of the roots of 
$$f(x):=\sum_{i=0}^n\omega_iV_{n-i}(K)(-x)^i,$$
then $0<a_1\leq r\leq R\leq a_n$, where $r$ and $R$ are the radii of $K$ relative to the $n$-dimensional unit ball $B_n$. We have confirmed this numerically for the orthoschemes of Theorem 2 up to dimension $n=21$. A further observation is that (in this case) the roots themselves appear to be real. We do not have an explanation for this.
\end{remark}

\section{Proof of Theorem 1 and a Conjecture}
The convex hull of a discretization of a crinkled arc $K_n=:\conv\{0, B(1/n), B(2/n),\ldots, 
B(1)\}$ is, up to scaling by $\sqrt{n}$, precisely the polytope $K$ in Theorem 2. Because 
\bqn
&&\lim_{n\rightarrow\infty}\left(\frac{1}{\sqrt{n}}\right)^{k/2}\sum_{l_1+l_2+\cdots+l_k\leq n}\frac{1}{\sqrt{l_1\,l_2\,\cdots l_k}}\\
&=&\int\!\!\int\!\!\cdots\!\!\int_{x_1+x_2+\cdots+x_k\leq 1}\frac{1}{\sqrt{x_1\, x_2\,\cdots x_k}}dx_1\;dx_2\cdots dx_k\\
&=&\omega_k \;,
\eqn
we have $V_k(K_B)\geq\omega_k/k!\,$. For the other direction of this inequality, it is enough to show it for any 
finite-dimensional convex body included in $K_B$. Because such a convex body can be approximated by convex polytopes included in $K_B$, and by the continuity of intrinsic volumes in finite-dimensions, it is enough to consider the case of a convex polytope $L_m$ with vertices $B(t_1)$, $B(t_2)$, $\dots$, $B(t_m)$, where the $t_i$'s are rational, say 
$t_i=k_i/n$, $1\leq i\leq m$. But then $L_m$ is a 
subset of $K_n$ as above, and this direction 
of the inequality is also clear.
\bigskip

We conclude with a conjecture. Following McMullen \cite{mcm2}, one knows that the sequence \linebreak
$m_k=(k+1)V_{k+1}/V_k,\,\,\,k=1,2,\ldots$ is decreasing in $k$. Recently it has been shown that limit 0 corresponds to the so--called $GC$ ($G$aussian $C$ontinuity) property: $V_1(K\cap B(t,\eps))\rightarrow 0$ as $\eps\rightarrow 0$ for every $t\in K$ (\cite{vit00}). $K_B$ is known to satisfy the latter condition, but only barely in a sense, and an easy consequence of Theorem 1 is that $m_k\sim{\rm constant}\,\cdot \, k^{-1/2}.$  We conjecture that this rate is extremal:

\bigskip
\noindent{\bf Conjecture} {\it For any infinite dimensional convex body with $V_1<\infty$, either 
$\lim_{k\rightarrow\infty}m_k>0$ or $m_k={\rm O}(k^{-1/2})$.}

\section{Acknowledgments} The authors thank R. Alexander for encouragement, L. Bobisud and S. Krone for valuable discussions, and L. Vitale for depictions and physical models. We also thank an anonymous referee for a close reading of the paper and for comments that improved the exposition.

\end{document}